\documentclass[review,authoryear]{elsarticle}

\usepackage{lineno,hyperref}

\usepackage[english]{babel}
\usepackage[utf8x]{inputenc}
\usepackage[T1]{fontenc}

\usepackage[a4paper,top=1cm,bottom=2cm,left=3cm,right=3cm,marginparwidth=1.75cm]{geometry}
\usepackage{indentfirst}
\usepackage{amsmath}
\usepackage{amssymb}
\usepackage{graphicx}
\usepackage[colorinlistoftodos]{todonotes}

\usepackage[linesnumbered,ruled,vlined]{algorithm2e}

\graphicspath{{images/}}

\modulolinenumbers[5]

\begin{document}

\title{Output maximization container loading problem with time availability constraints}





\author[aus,br]{Pedro B. Castellucci\corref{mycorrespondingauthor}}
\cortext[mycorrespondingauthor]{Pedro B. Castellucci}
\ead{pbc@icmc.usp.br}

\author[br]{Franklina M. B. Toledo}
\ead{fran@icmc.usp.br}

\author[aus]{Alysson M. Costa}
\ead{alysson.costa@unimelb.edu.au}

\address[aus]{The University of Melbourne}
\address[br]{Universidade de São Paulo}

\begin{frontmatter}
\begin{abstract}
Research on container loading problems has been proved effective in increasing the filling rate of containers in different practical situations. However, the broader logistic context might pressure the loading process, leading to sub-optimal solutions. Some facilities like cross-docks have reduced storage space which might force early loading activities. We propose a container loading problem which accounts for this limited storage by explicitly considering the schedule of arrival for the boxes and the departure time of the trucks. Also, we design a framework which handles the geometric and temporal characteristics of the problem separately, enabling the use of methods found in the literature for solving the extended problem. Our framework can handle uncertainty in the schedule and be used to quantify the impact of delays on capacity utilization and departure time of trucks.
\end{abstract}


\end{frontmatter}

\section{Introduction}
\label{sec:intro}

Container loading problems have the goal of improving the effectiveness of logistic systems by increasing the capacity utilization of the trucks/containers. Its most popular version consists in selecting a subset of items (boxes) that maximize the volume (or value) loaded into a limited number of containers. However, the logistic environment can pressure the loading process to the use of less effective solutions. That is the case with cross-docking. Cross-docks are facilities with limited storage space which operate by synchronizing the inbound and outbound flow of trucks to potentially improve rates of consolidation, delivery and lead times and costs \citep{VanBelle2012}. 

One of the most popular optimization problems in the cross-docking literature is the scheduling of inbound and outbound trucks at its docks. To provide insights on promising solution methodologies and behavior of the system, some researchers have focused on cases with one inbound and one outbound dock \citep{Yu2008,Chen2009,Vahdani2010,BolooriArabani2010,BolooriArabani2011,Liao2012,Amini2016, Keshtzari2016}, but cases with multiple inbound and outbound docks have been becoming more frequent \citep{Bellanger2013, VanBelle2013, ladier2014crossdock,Assadi2016,Cota2016, Wisittipanich2017}. The scheduling problem is usually modeled with the goal of minimizing total operation time \citep{Yu2008,Chen2009,Vahdani2010,Liao2012,Keshtzari2016,Bellanger2013,Cota2016,Wisittipanich2017,BolooriArabani2011meta,Shakeri2012,Gelareh2016,Bazgosha2017}. It is also possible to find examples of other metrics: Storage (volume or cost) or material handling \citep{ladier2014crossdock,RahmanzadehTootkaleh2016,Maknoon2017}, travel distance \citep{Chmielewski2009} and earliness and/or tardiness alone \citep{BolooriArabani2010,Assadi2016} or combined with others such as travel time \citep{VanBelle2013} and probability of breakdowns \citep{Amini2016}.

The combination of limited storage space and pressure for the early dispatching of trucks in a cross-dock imposes additional challenges to the resulting container loading problems. The container loading literature already contemplates a wide range of practical constraints related to the weight of the containers \citep{davies1999weight,Ceschia2013,Lim2013}, the position of individual boxes \citep{Bortfeldt2012,Zheng2014,Wei2015,Jamrus2016,Sheng2016,Tian2016} and the arrangement of the boxes \citep{bischoff1995issues,Junqueira2012a,Takahara2005,alonso2018}.  For an extensive review on practical constraints and solution methodologies see \citep{Bortfeldt2013} and \citep{Zhao2016}, respectively. Although these constraints might be present in container loading in cross-docking environments, we are not aware of any research dealing with the specific needs of cross-docks: (i) considering the loading time of trucks and (ii) time availability of the boxes, due to limited storage.

We explore these two specific challenges of container loading problem in cross-docking environments: we consider the problem of selecting a subset of boxes to maximize the loaded volume and minimize the ready time of the truck/container. Furthermore, we mind the schedule of arrival of boxes. Since the availability of the boxes in time might be uncertain, we propose a stochastic dynamic programming framework which can accommodate this uncertainty and take advantage of the existing methods to solve the geometric aspects of the problem. 

\section{Container loading problem with time availability constraints (CLPTAC)}
\label{sec:problem-def}

We propose an output maximization version of a container loading problem with time availability constraints (CLPTAC-Om). The objective is to position a subset of boxes inside a containers maximizing the loaded volume and minimizing the ready time of the container to be dispatched. The boxes must not overlap with each other, must be fully inside the container and cannot be loaded before they are available. Also, we assume boxes and containers' faces must be parallel to each other (orthogonal packing). We define deterministic and stochastic versions of the problem with respect to the arrival time of the boxes. Moreover, we assume that boxes leaving the dock go to the same destination, which is a common situation in cross-docks with destination exclusive dock holding patterns (\cite{Ladier2016review}).

More formally, assume a finite time horizon $[0, T]$. Let $\mathcal{B}$ be the set of boxes, each box $i \in \mathcal{B}$ is defined by its length along each axis $(l_{i1}, l_{i2}, l_{i3})$ and the time it becomes available for loading $t_i$. Simulating the arrival process of the boxes, we assume that, once a box becomes available, it remains available until it is shopped to its destination. We consider a single container with dimensions $(L_1, L_2, L_3)$. For the deterministic version, the arrival time of boxes $t_i \in [0, T]$, $i \in \mathcal{B}$, is defined \textit{a priori}. For the stochastic version, only the probability distribution of $t_i$ is known.The remainder of this paper focus on the stochastic version, although the dynamic programming framework we propose (Section \ref{sec:stoch-dyn-prog-alg}) can also be applied to the deterministic version.

\section{A dynamic programming framework for CLPTAC-Om}
\label{sec:stoch-dyn-prog-alg}

The framework we propose to solve the stochastic CLPTAC-Om decomposes the temporal and geometric aspects of the problem. This enables the use of any approach for the classical container loading problems to compute the costs of the states in the dynamic programming algorithm. For this, we uniformly discretize the time horizon $[0, T]$ into a set of time periods $\mathcal{T}$ and define $\mathcal{B}_t$, $t \in \mathcal{T}$, to be the set of boxes available at time period $t \in \mathcal{T}$. An effective discretization might not be trivial for general systems, but for cross-docks a simple uniform discretization is justifiable due to usual high and stable demands.

 The dynamic programming approach is formulated as an optimum stopping problem \citep{peskir2006optimal}. We want to decide when to stop waiting for new shipments ($\mathcal{B}_t,\ t \in \mathcal{T}$) and load the boxes into the truck. The trade-off is: the longer we wait, the more boxes we have available for packing yielding a potentially  better volume occupation; on the other hand, the sooner we stop waiting (and perform the loading) the smaller is the ready time for the selected boxes.

Let the state of the system be the set of boxes $\mathcal{X}_t$ available at period $t \in \mathcal{T}$ and   $\mathcal{C}(\mathcal{X}_t, s_t)$ be the cost of decision $s_t \in \{0, 1\}$, $t \in \mathcal{T}$, in state $\mathcal{X}_t$. The possible decisions are to perform the loading or to wait for more boxes, henceforth coded as 0 (load) and 1 (wait). Also, let $\mathcal{E}(\mathcal{X}_t)$ be a solution of a traditional container loading problem with boxes in $\mathcal{X}_t$, namely, a solution that maximizes the volume loaded with boxes represented by $\mathcal{X}_t$ without any temporal constraints and $\overline{\mathcal{E}} (\mathcal{X}_t)$ be the volume of boxes that were not loaded. Finally, let $\mu$ be a parameter indicating the cost of empty volume in the container. Thereby, with $\mathcal{T}^*$ being the set of all time periods but the last, the cost functions of the dynamic programming algorithm are given by (\ref{eq:cost-u0})--(\ref{eq:finish-process}). 

 \begin{align}
 \label{eq:cost-u0}
 & \mathcal{C}(\mathcal{X}_t, 0) = \frac{1}{|\mathcal{T}|} + \mu  \Bigg(1 - \frac{\mathcal{E}(\mathcal{X}_t)}{L_1L_2L_3}\Bigg), \quad  && t \in \mathcal{T}^*, \mathcal{X}_t \neq \emptyset,\\
 \label{eq:cost-u1}
 & \mathcal{C}(\mathcal{X}_t, 1) = \frac{1}{|\mathcal{T}|}, && t \in \mathcal{T}^*, \mathcal{X}_t \neq \emptyset,\\
 \label{eq:cost-uN}
 & \mathcal{C}(\mathcal{X}_{|\mathcal{T}|}) = \frac{1}{|\mathcal{T}|} + L_1L_2L_3 \Bigg \lceil \frac{\overline{\mathcal{E}}(\mathcal{X}_{|\mathcal{T}|})}{L_1L_2L_3}\Bigg \rceil, &&\\
 \label{eq:finish-process}
 & \mathcal{C}(\emptyset, s) = \mathcal{C}(\emptyset) = 0,  && t \in \mathcal{T}^*, s \in \{0, 1\}.
\end{align}
Equations (\ref{eq:cost-u0}) define the cost of the decision to stop the waiting process and proceed with the loading with available boxes $\mathcal{X}_t$ at time $t \in \mathcal{T}^*$. If the decision is to wait for more boxes, the respective partial cost is given by (\ref{eq:cost-u1}). The terminal cost, i.e., the cost of reaching the last period in the planning horizon is defined by a bound on the volume of the trucks needed to transport all remaining boxes (\ref{eq:cost-uN}), in a sense, it is the inventory holding costs, which should be high to preserve the cross-docking philosophy of reduced inventory. If there are no boxes available for loading, $\mathcal{X}_t = \emptyset$, the costs are given by~(\ref{eq:finish-process}). We want to optimize the objective (\ref{eq:clptac-input-dyn-obj}), in which $\mathbb{E}$ denotes the expected value due to the uncertainty in the arrival times of the boxes.

\begin{equation}
  \label{eq:clptac-input-dyn-obj}
  \mbox{Min }\mathbb{E}\Bigg(\sum_{t \in \mathcal{T}^*}( C(\mathcal{X}_t, s_{t})) + C(\mathcal{T})\Bigg).
\end{equation}
Which is equivalent to minimizing the sum of the ready times (in time periods) and the cost of the empty volume that it is possible to achieve with the boxes available until we decide to stop waiting, $s_t = 0$. Note that, if the decision in a particular period is to wait, $s_t = 1$, we sum the cost (regarding ready time) of that period, $\frac{1}{|\mathcal{T}|}$. Otherwise, $s_t = 0$, we still sum the cost related to the ready time (plus the cost of empty volume). In other words, we assume that the loading process only ends at the end of a period. We use $|\mathcal{T}|$ and $L_1L_2L_3$ as normalizing factor for ready time and volume, respectively. 

Then, we can find an optimal solution for the stochastic case by using the value functions (\ref{eq:dyn-prob-recursion-stochastic}) and(\ref{eq:dyn-prob-recursion-end-stochastic}).

\begin{align}
\label{eq:dyn-prob-recursion-stochastic}
& \mathcal{V}_t(\mathcal{X}_t)= \mbox{Min} \Bigg\{ \frac{1}{|\mathcal{T}|} + \sum_{\overline{b} \subseteq \mathcal{X}_{t+1} }  \rho_{t+1}( \overline{b}\ |\ \mathcal{X}_t) \ \mathcal{V}_{t+1}(\overline{b}),\ \frac{1}{|\mathcal{T}|} + \mu \Bigg(1 - \frac{\mathcal{E}(\mathcal{X}_t)}{L_1L_2L_3}\Bigg) \Bigg\},  t \in \mathcal{T}^*, \mathcal{X}_t \subseteq \overline{\mathcal{B}}_t,\\
\label{eq:dyn-prob-recursion-end-stochastic}
& \mathcal{V}_{|\mathcal{T}|}(\mathcal{X}_{t}) = \frac{1}{|\mathcal{T}|} +L_1L_2L_3 \Bigg \lceil \frac{\overline{\mathcal{E}}(\mathcal{X}_{|\mathcal{T}|})}{L_1L_2L_3}\Bigg \rceil, \quad \mathcal{X}_t \subseteq \overline{\mathcal{B}}_{|\mathcal{T}|},
\end{align}
in which $\rho_t (\mathcal{B}_t)$ is the probability of shipment $\mathcal{B}_t$ arriving in time period $t$ and $\overline{\mathcal{B}}_t$ is the set of all shipments $\{\mathcal{B}_1, \ldots, \mathcal{B}_t\}$ that could be available at time $t$. Thus,
the expected value of $\mathcal{V}_1(\mathcal{X}_1)$ gives an optimum solution to the stochastic version of the problem.

To compute $\rho_{t+1}(\overline{b}\ | \ \mathcal{X}_t)$ in (\ref{eq:dyn-prob-recursion-stochastic}), we need a probabilistic model for the arrival times of boxes in $\mathcal{B}_t$. To illustrate our approach, we assumed the arrival time of subset $\mathcal{B}_{t_1}$ is independent of the arrival time of $\mathcal{B}_{t_2}$, $t_1 \neq t_2$, and the probability of a subset $\mathcal{B}_t$ arriving at time $i \in \mathcal{T}$ is given by:

\begin{equation}
\label{eq:random-definition}
\rho_i(\mathcal{X}_t = \mathcal{B}_t) = 
	\begin{cases} 
    	\alpha_{it} (1 - \alpha_{it})^{i-t}, &\mbox{if } i \geq t, \\ 
		0 & \mbox{if } i < t,
	\end{cases}
\end{equation}
in which $\alpha_{it} \in (0, 1]$, $i, t \in \mathcal{T}$, is a probability index associated with subset $\mathcal{B}_t$ arriving at time $i \in \mathcal{T}$, $i \geq t$. Within this framework, $\alpha_{it}$, can be interpreted as a reliability profile for the company in charge of delivering the subset (cargo) $\mathcal{B}_t$. If $\alpha_{it} = 1$, $i,t \in \mathcal{T}$, we have the deterministic case, in which we known the arrival time of each box \textit{a priori}. For simplicity, we used $\alpha_{it} = \alpha$, $i \in \mathcal{B}_t$, $t \in \mathcal{T}$, which can be viewed as a reliability index for the corresponding logistic environment. 

\begin{figure}[htpb]
\centering
\includegraphics[scale=0.65,trim={0 1.5cm 0 0.5cm}]{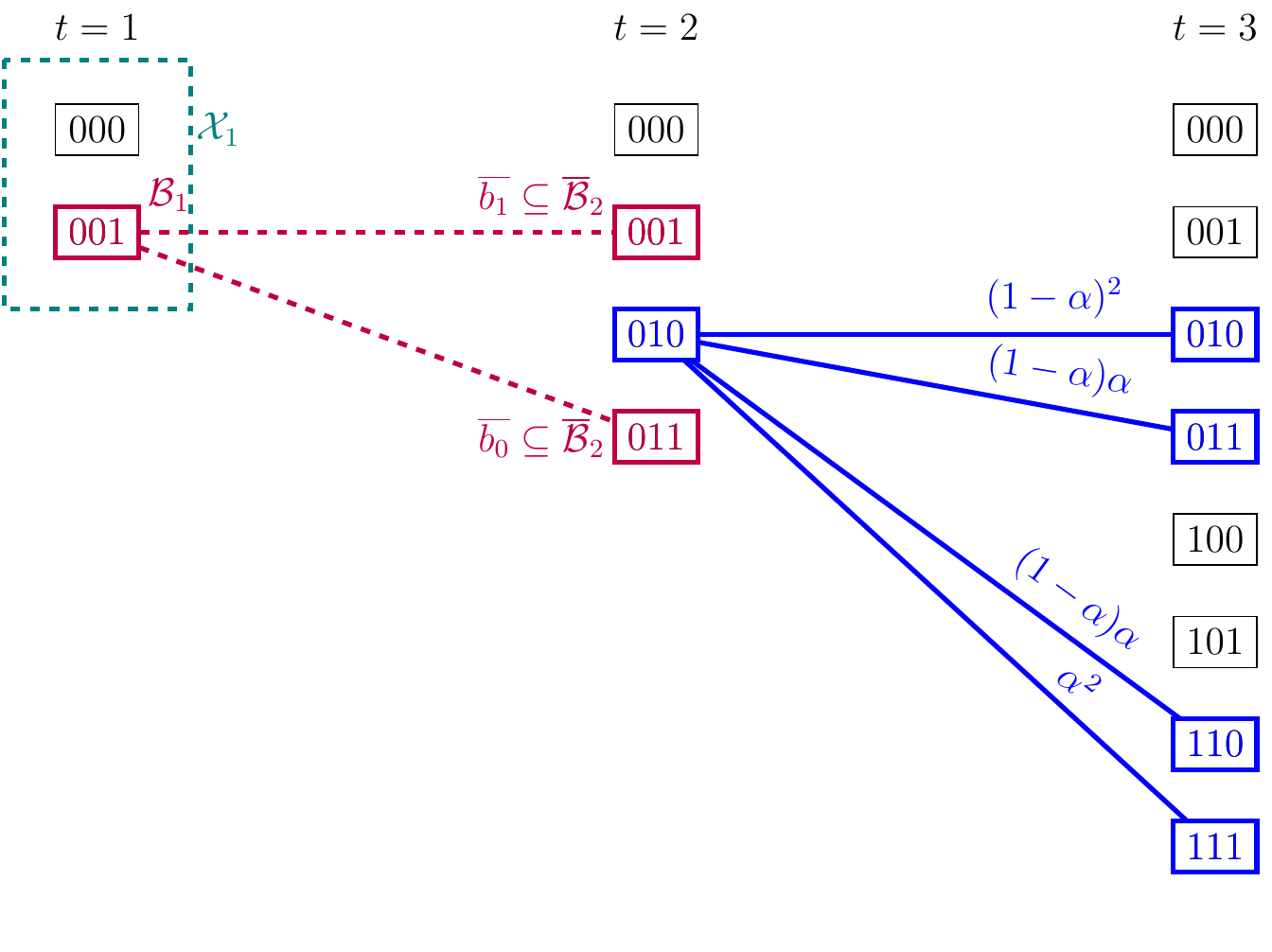}
\caption{Example of a scenario with $\mathcal{T}= \{1, 2, 3\}$. It is a graphical example of equation (\ref{eq:dyn-prob-recursion-stochastic}). The dashed red lines show possible next states for $\mathcal{X}_1$. The full blue lines show the states which need to be evaluated at state $010$ at time $t=2$.}
\label{fig:dynProgImpl}
\end{figure}

Note that equations (\ref{eq:dyn-prob-recursion-stochastic}) and (\ref{eq:dyn-prob-recursion-end-stochastic}) are independent of the probabilistic model for the arrival of the boxes.
Figure \ref{fig:dynProgImpl} exemplifies a scenario with $\mathcal{T} = \{1, 2, 3\}$. Each cell is a bit array ($\mathcal{B}_3\mathcal{B}_2\mathcal{B}_1$) representing the current state, i.e., set of boxes that have actually arrived at a point in time. The bit array $\overline{b}_1 = 001$, for example, means that shipment 1 has arrived while shipments 2 and 3 have not. In each time period, for each bit array with non-zero probability of occurrence, we need to evaluate the decisions of waiting and loading (equations (\ref{eq:dyn-prob-recursion-stochastic})). The full blue lines represent the computation of the sum in equation (\ref{eq:dyn-prob-recursion-stochastic}) for all $\rho_i (\mathcal{X}_3 = \mathcal{B}_3) \neq 0$.

\subsection{Computational experiments}
\label{sec:clptac-output-comp-exper}

The benchmark problems for the CLPTAC-Om were generated based on the 47 instances proposed by \cite{Ivancic1989} for the traditional container loading problem. In this benchmark, instances have from two to five types of boxes, the number of boxes varies from 47 to 181 and we consider only one container per instance -- containers have different dimensions across instances. For each box, we randomly generated a time from which it becomes available following a discrete uniform probability distribution in $[1, |\mathcal{T}|], |\mathcal{T}| = 10$. Therefore, we are assuming a planning horizon with ten time slots (equivalent to 40 minutes slots in an eight-hour working day).

We experimented with different trade-off parameters for the cost of empty spaces $\mu \in \{1, 2, 4\}$. To solve the sub-problems, we used the container loading model defined in \cite{Chen1995} without allowing the rotation of the boxes. We solved this model with Gurobi 7.0 in Ubuntu 12.04 using four cores of an Intel\textsuperscript{\textregistered} Xeon CPU E5-2620 at a frequency of 2 GHz each. In the following, we denote CLPTAC$\mu$ instances in which the cost of empty spaces is $\mu \in \{1, 2, 4\}$.

\subsection{Results on the dynamic programming approach}
\label{sec:clptac-res-output-dynprog}

Solving the stochastic model for an instance of the CLPTAC-Om may be computationally expensive. Each evaluation of function $\mathcal{V}_t (\cdot)$ needs to solve a deterministic version of the container loading problem (via $\mathcal{E}(\cdot)$), which for real-world cases can be impractical itself, using an exact method. Therefore, to validate our framework, we solved the sub-problem in which all items are available ($\mathcal{X}_{|\mathcal{T}|} = \cup_{t \in \mathcal{T}}\mathcal{B}_{t}$) (with a time limit of 3600 seconds) and then reduced the time limit to five seconds to find a solution for sub-problems related to $V_t$, $t \in \mathcal{T}^*$ (we noted that the solver is able to find ``good'' solutions quickly and then struggles to prove optimality). The goal of these experiments is to evaluate the effect of the cost of empty volume ($\mu$) and of the reliability index ($\alpha$) on total ready time and volume occupation. Then, with the optimal policy from equation (\ref{eq:dyn-prob-recursion-stochastic}), we simulated the arrival of boxes 5000 times for each $\alpha \in \{0.5, 0.6 \ldots, 0.9\}$.

\begin{figure}[htpb]
\centering
\includegraphics[scale=0.32]{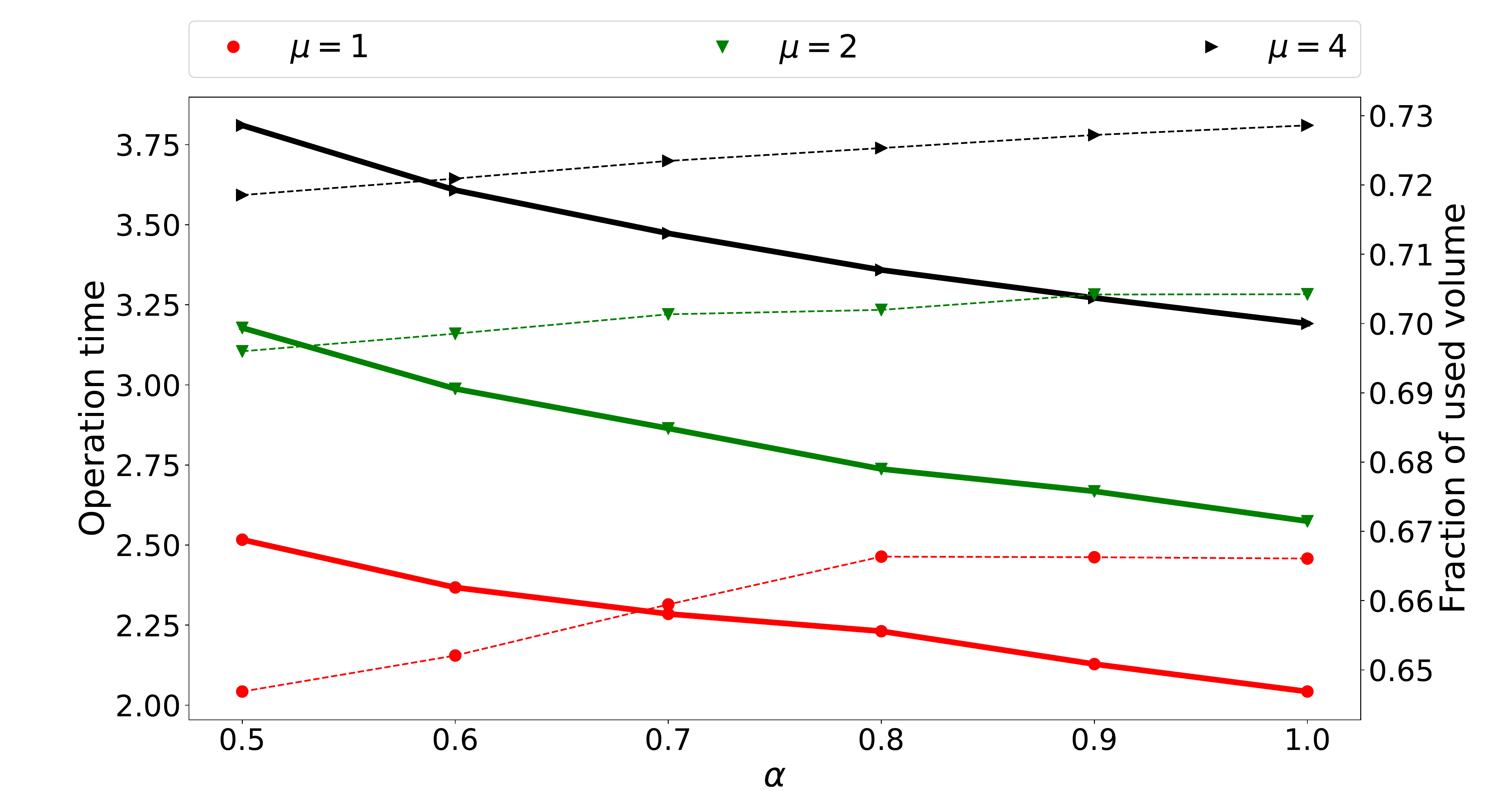}
\caption{Profile of loaded volume (dashed line) and operation time (full line) related to the reliability index ($\alpha$) for $\mu \in \{1, 2, 4\}$. The values were averaged across all instances.}
\label{fig:res-stochastic-dp-time-volume}
\end{figure}

The experiments revealed that operation time and volume occupation tend to improve as reliability index increases. Although it was not possible to note a significant statistical difference among cases with different $\mu$, an increasing trend can be observed in both ready time and volume occupation (Figure \ref{fig:res-stochastic-dp-time-volume}). Regarding operation time, this trend can be explained due to a higher risk of waiting in a environment with lower reliability index. The same reasoning applies to the volume occupation, the higher the risk of waiting the lower the occupation achieved. For the decision marker, it suggests that the more reliable the logistic environment the better the expected vehicles occupation and the operation time.

Also, the experiments revealed a linear trend specially for the operation time in relation to the reliability index. It can help managers in devising business plans to account for costs of operation, delivery times and value of possible penalties for delay deliveries in seasonal variations of the reliability index. This, besides allowing a more efficient management of the facility, can potentially increase the service level for final customers.

\section{Conclusions and future research}
\label{sec:clptac-output-conclusions}

High filling rates of trucks/container are essential for effective road distribution systems. However, achieving these rates can be challenging due to the pressure for lower and lower delivery times. We define an \textit{Output Maximization Container Loading Problem with Time Availability Constraints} which accommodates these two opposing goals: high filling rates and low delivery times. Moreover, we propose a dynamic programming framework which is suited to handle the problem, including uncertainty on the availability of boxes in time. This framework separates the geometric and temporal aspects of the problem, enabling the use of all available solution methodologies for tackling different practical constraints in the container loading process. Future work may extend the core decision problem considered here in order to include additional decisions in the context of cross-docking strategies. In particular, some natural extensions would be the integration of scheduling decisions for the arrival of shipments and  the consideration of multiple docks. Other avenue for research is the the consideration of multiple destinations, with the definition of shipping routes in the second layer of the distribution network.

\section*{Acknowledgements}

We thank State of São Paulo Research Foundation (FAPESP), process numbers, \#2017/01097-9, \#2015/15024-8, CEPID \#2013/07375-0 and National Council for Scientific and Technological Development (CNPq), grant \#306918/2014-5 from Brazil.

We also thank Professor Eduardo Fontoura Costa from the Instituto de Ciências Matemática e de Computação (ICMC) of the University of Sao Paulo for the initial discussions regarding the dynamic programming algorithm.

\section*{References}
\bibliography{bib/all}
\bibliographystyle{plainnat}

\end{document}